\documentclass[12pt]{article}

\usepackage{amsmath}
\usepackage{amssymb}
\usepackage{amsthm}
\usepackage[pdfpagemode=None]{hyperref}
\usepackage{graphicx}
\usepackage{pstricks}

\newtheorem{thm}{Theorem}[section]

\newtheorem{lem}[thm]{Lemma}

\newtheorem{assumption}[thm]{Assumption}

\newtheorem{definition}[thm]{Definition}

\newtheorem{example}[thm]{Example}
\newenvironment{exmp}{\begin{example}\rm}{\end{example}}
\newtheorem{remark}[thm]{Remark}
\newenvironment{rem}{\begin{remark}\rm}{\end{remark}}
\newtheorem{tab}{Table}

\def\eps{\varepsilon}

\title{A note on a complex Hilbert metric with application to domain of
analyticity for entropy rate of hidden Markov processes}

\author{{\small 
\begin{tabular}{ccc}
Guangyue Han& Brian Marcus & Yuval Peres\\
Department of Mathematics&Department of Mathematics&Department of Statistics\\
University of Hong Kong&University of British Columbia&University of California, Berkeley\\
Email:ghan@maths.hku.hk&Email:marcus@math.ubc.ca&Email: peres@stat.berkeley.edu\\
\end{tabular}}}

\date{{\normalsize \today}}

\begin{document}\maketitle\thispagestyle{empty}


\vspace{1cm}

\centerline{DRAFT}

\begin{abstract}
In this note, we show that small complex perturbations of positive
matrices are contractions, with respect to a complex version of the
Hilbert metric, on the standard complex simplex. We show that this
metric can be used to obtain estimates of the domain of analyticity
of entropy rate for a hidden Markov process when the underlying
Markov chain has strictly positive transition probabilities.
\end{abstract}

The purpose of this note is twofold.  First, in
Section~\ref{complex-Hilbert}, we introduce a complex version of the
Hilbert metric on the standard real simplex. This metric is defined
on a complex neighbourhood of the interior of the standard real
simplex, within the standard complex simplex.  We show that if the
neighbourhood is sufficiently small, then for any sufficiently small
complex perturbation of a strictly positive square matrix acts as a
contraction, with respect to this metric. While this paper was
nearing completion, we were informed of a different complex Hilbert
metric, which was recently introduced. We briefly discuss the
relation between this metric~\cite{Dubois} and our metric in
Remark~\ref{Dubois}.

Secondly, we show how one can use a complex Hilbert metric to obtain
lower estimates of the domain of analyticity of entropy rate for a
hidden Markov process when the underlying Markov chain has strictly
positive transition probabilities.  The domain of analyticity is
important because it specifies an explicit region where a Taylor
series converges to the entropy rate and also gives an explicit
estimate on the rate of convergence of the Taylor approximation.

In principle, an estimate on the domain can be obtained by examining
the proof of analyticity in~\cite{gm05}. That proof was based on a
contraction mapping argument, using the fact that the real Euclidean
metric is equivalent to the real Hilbert metric. In
Section~\ref{background}, we revisit certain aspects of the proof
and outline how to modify the proof using a complex Hilbert metric;
this yields a more direct estimate. In Section~\ref{example}, we
illustrate this with a small example, using our Hilbert metric.

We remark that the entropy rate of a hidden Markov process can be
interpreted as a top Lyapunov exponent for a random matrix product
~\cite{Holliday}. In principle, a complex Hilbert metric can be
used, more generally, to estimate the domain of analyticity of the
top Lyapunov exponent for certain random matrix products; see
~\cite{Peres1},~\cite{Peres2}.

\section{Complex Hilbert Metric}
\label{complex-Hilbert}

We begin with a review of the real Hilbert metric. Let $B$ be a
positive integer, and let $W$ be the standard simplex in
$B$-dimensional real Euclidean space:
$$
W = \{w=(w_1, w_2, \cdots, w_B) \in \mathbb{R}^B:w_i \geq 0, \sum_i
w_i=1\},
$$
and let $W^{\circ}$ denote its interior, consisting of the vectors
with positive coordinates. For any two vectors $v, w \in W^{\circ}$, the Hilbert
metric~\cite{se80} is defined as
\begin{equation}  \label{HilbertMetric}
d_H(w, v)=\max_{i, j} \log \left( \frac{w_i/w_j}{v_i/v_j} \right).
\end{equation}
For a $B \times B$ strictly positive matrix $T=(t_{ij})$, the
mapping $f_T$ induced by $T$ on $W$ is defined by $f_T(w) = wT/ (w T
{\bf 1} )$, where $\bf 1$ is the all-ones vector.  It is well known
that $f_T$ is a contraction mapping under the Hilbert
metric~\cite{se80}. The contraction coefficient of $T$, which is
also called the Birkhoff coefficient, is given by:
\begin{equation}  \label{BirkhoffCoefficient}
\tau (T)=\sup_{v \neq w} \frac{d_H(vT, wT)}{d_H(v, w)}=\frac{1-\sqrt{\phi(T)}}{1+\sqrt{\phi(T)}},
\end{equation}
where $\phi(T)=\min_{i, j, k, l} \frac{t_{ik}t_{jl}}{t_{jk}t_{il}}$.
This result extends to the case where $T$ has all columns strictly positive or all zero
and at least one strictly positive column (then, in the definition of $\phi(T)$,
consider only $k,l$ corresponding to strictly positive columns).

Let $W_{\mathbb{C}}$ denote the complex version of $W$, i.e.,
$W_{\mathbb{C}}$ denotes the complex simplex comprising the vectors
$$
\{w=(w_1, w_2, \cdots, w_B) \in \mathbb{C}^B : \sum_i w_i=1\}.
$$

Let $W_{\mathbb{C}}^+ = \{v \in W_{\mathbb{C}}: \mathcal{R}(v_i) >0\}$.  For
$v,w \in W_{\mathbb{C}}^+$, let
\begin{equation}  \label{ComplexHilbertMetric} d_H(v, w)=\max_{i
, j} \left| \log \left( \frac{w_i/w_j}{v_i/v_j} \right) \right|,
\end{equation}
where $\log$ is taken as the principal branch of the complex
$\log(\cdot)$ function (i.e., the branch whose branch cut is the
negative real axis).  Since the principal branch of $\log$ is
additive on the right-half plane, $d_H$ is a metric on
$W_{\mathbb{C}}^+$, which we call a {\em complex Hilbert metric}.

We begin with the following very simple lemma.

\begin{lem}  \label{maxone}
Let $n \geq 2$. For any fixed $z_1, z_2, \cdots, z_n, z \in
\mathbb{C}$ and fixed $t > 0$,
we have
$$
\sup_{t_1, \cdots, t_n \geq 0,~ t_1+t_2+\cdots+t_n=t} ~\left| t_1
z_1 + t_2 z_2+ \cdots+t_n z_n+z \right| = \max_{i=1,\cdots, n}
\left| t z_i + z \right|.
$$
\end{lem}

\begin{proof}

The convex hull of $z_1, z_2, \cdots, z_n$ is a solid polygon,
taking the form
$$
\{(t_1/t) z_1 + (t_2/t) z_2+ \cdots+(t_n/t) z_n: t_1, t_2, \cdots,
t_n \geq 0, t_1+t_2+\cdots+t_n=t\}.
$$
By convexity, the distance from any point in this solid polygon to
the point $(-1/t)z$ will achieve the maximum at one of the extreme
points, namely
$$
\sup_{t_1, \cdots, t_n \geq 0, t_1+t_2+\cdots+t_n=t} \left| (t_1/t) z_1 + (t_2/t) z_2+
\cdots+(t_n/t) z_n-(-1/t) z \right| = \max_{i=1,\cdots, n} \left|
z_i - (-1/t) z \right|.
$$
The lemma then immediately follows.

\end{proof}

The following lemma is implied by the proof of Lemma 2.1
of~\cite{ro03}; we give a proof for completeness.

\begin{lem}  \label{aA}
For fixed $a_1, a_2, \cdots, a_B > 0 \in \mathbb{R}$ and fixed $x_1, x_2 \cdots, x_B > 0 \in \mathbb{R}$,
define:
$$
D_n=\frac{a_n x_n}{\sum_{m=1}^B a_m x_m}-\frac{x_n}{\sum_{m=1}^B x_m}.
$$
Let $\mathcal{T}_0=\{n: D_n \geq 0\}$ and $\mathcal{T}_1=\{n: D_n < 0\}$. Then we have
$$
\sum_{n \in \mathcal{T}_0} D_n=\sum_{n \in \mathcal{T}_1} |D_n| \leq \frac{1-\sqrt{a/A}}{1+\sqrt{a/A}},
$$
where $a=\min \{a_1, a_2, \cdots, a_B\}$ and $A=\max \{a_1, a_2, \cdots, a_B\}$.
\end{lem}

\begin{proof}
It immediately follows from $\sum_{n=1}^B D_n=0$ and the definitions of $\mathcal{T}_0$ and $\mathcal{T}_1$ that
$$
\sum_{n \in \mathcal{T}_0} D_n=\sum_{n \in \mathcal{T}_1} |D_n|.
$$
Now
$$
\sum_{n \in \mathcal{T}_0} D_n = \sum_{n \in \mathcal{T}_0} \left(
\frac{ a_n x_n}{\sum_{m \in \mathcal{T}_0} a_m x_m+\sum_{m \in
\mathcal{T}_1} a_m x_m}-\frac{ x_n}{\sum_{m \in \mathcal{T}_0}
x_m+\sum_{m \in \mathcal{T}_1} x_m} \right)
$$
$$
\leq \sum_{n \in \mathcal{T}_0} \left( \frac{A x_n} {A \sum_{m \in
\mathcal{T}_0} x_m+a \sum_{m \in \mathcal{T}_1}
x_m}-\frac{x_n}{\sum_{m \in \mathcal{T}_0} x_m+\sum_{m \in
\mathcal{T}_1} x_m} \right)
$$
Let
$$
z=\frac{\sum_{m \in \mathcal{T}_1} x_m}{\sum_{m \in \mathcal{T}_0}
x_m},
$$
we then have
$$
\sum_{n \in \mathcal{T}_0} D_n \leq \frac{1}{1+(a/A) z}-\frac{1}{1+z}=f(z).
$$
Simple calculus shows that $f(z)$ will be bounded above by
$\frac{1-\sqrt{a/A}}{1+\sqrt{a/A}}$ on $[0, \infty)$. This
establishes the lemma.
\end{proof}

Let $W_{\mathbb{C}}^{\circ}(\delta)$ denote the ``relative''
$\delta$-neighborhood of $W^{\circ}$ in $W_{\mathbb{C}}$, i.e.,
$$
W_{\mathbb{C}}^{\circ}(\delta)=\{v=(v_1, v_2, \cdots, v_B) \in
W_{\mathbb{C}}: \exists u \in W^{\circ}, |v_i-u_i| \leq \delta
|u_i|, i=1, 2, \cdots, B\}.
$$
Note that when $\delta <1$,  $W_{\mathbb{C}}^{\circ}(\delta) \subset
W_{\mathbb{C}}^+$ and so the complex Hilbert metric is defined on
$W_{\mathbb{C}}^{\circ}(\delta)$.

We consider complex matrices $\hat{T}=(\hat{t}_{ij})$ which are
perturbations of a strictly positive matrix $T=(t_{ij})$. For such a matrix $T$ and $r > 0$, let $B_T(r)$ denote the set of all complex matrices
$\hat{T}$ such that for all $i, j$,
$$
|t_{ij}-\hat{t}_{ij}| \leq r.
$$

With the aid of the above lemmas, we shall prove:
\begin{thm}   \label{ComplexContraction}
Let $T$ be a strictly positive matrix.  There exist $r, \delta > 0$
such that whenever $\hat{T} \in B_T(r)$, $f_{\hat{T}}$ is a
contraction mapping on $W_{\mathbb{C}}^{\circ}(\delta)$ under the
complex Hilbert metric.
\end{thm}

\begin{proof}
For $\hat{x}, \hat{y} \in W_{\mathbb{C}}^+$, $\hat{x} \ne \hat{y}$,
and $i,j$, let
$$ L_{ij} = \frac{\log(\sum_m \hat{x}_m \hat{T}_{mi}/\sum_m \hat{x}_m
\hat{T}_{mj})-\log(\sum_m \hat{y}_m \hat{T}_{mi}/\hat{y}_m
\hat{T}_{mj})}{\max_{k, l} \left| \log (\hat{x}_k/\hat{y}_k) -
\log(\hat{x}_l/\hat{y}_l) \right|}.
$$
Note that $$ \frac{d_H(\hat{x} \hat{T}, \hat{y}
\hat{T})}{d_H(\hat{x}, \hat{y})}= \max_{i,j} |L_{ij}|.
$$
It suffices to prove that there exists $0 < \rho <1$ such that for sufficiently small $r, \delta
> 0$, $\hat{x}, \hat{y} \in W_{\mathbb{C}}^{\circ }(\delta)$,
$\hat{x} \ne \hat{y}$, $\hat{T} \in B_T(r)$, and any $i, j$,
$$
|L_{ij}| < \rho.
$$
For each $m$, let $\hat{c}_m=\log \hat{x}_m/\hat{y}_m$; then
$\hat{x}_m=\hat{y}_m e^{\hat{c}_m}$. Choose $p \neq q$ such that
$$
\left| \hat{c}_p - \hat{c}_q \right|=\max_{k, l} \left| \hat{c}_k - \hat{c}_l \right|.
$$
Hence:
$$
L_{ij} = \frac{\log(\sum_m \hat{y}_m e^{\hat{c}_m-\hat{c}_q}
\hat{T}_{mi}/\sum_m \hat{y}_m e^{\hat{c}_m-\hat{c}_q}
\hat{T}_{mj})-\log(\sum_m \hat{y}_m \hat{T}_{mi}/\hat{y}_m
\hat{T}_{mj})}{|\hat{c}_p - \hat{c}_q|}.
$$
Define
$$
F(t)=\log(\sum_m \hat{y}_m e^{(\hat{c}_m-\hat{c}_q)t}
\hat{T}_{mi}/\sum_m \hat{y}_m e^{(\hat{c}_m-\hat{c}_q)t}
\hat{T}_{mj}).
$$
Since
$$
|F(1)-F(0)|=\left| \int_0^1 F'(t) dt \right| \leq \max_{\xi \in [0,
1]} |F'(\xi)|,
$$
we have
\begin{equation}
\label{Lijbound}
|L_{ij}|=\frac{|F(1)-F(0)|}{|\hat{c}_p-\hat{c}_q|} \leq
\frac{\max_{\xi \in [0, 1]} |F'(\xi)|}{|\hat{c}_p-\hat{c}_q|}.
\end{equation}
Note that $F'(\xi)$ takes the following form:
$$
F'(\xi)=\frac{\sum_m (\hat{c}_m-\hat{c}_q)  \hat{y}_m
e^{(\hat{c}_m-\hat{c}_q) \xi} \hat{T}_{mi}}{\sum_m \hat{y}_m
e^{(\hat{c}_m-\hat{c}_q) \xi} \hat{T}_{mi}}-\frac{\sum_m
(\hat{c}_m-\hat{c}_q) \hat{y}_m e^{(\hat{c}_m-\hat{c}_q) \xi}
\hat{T}_{mj}}{\sum_m \hat{y}_m e^{(\hat{c}_m-\hat{c}_q) \xi}
\hat{T}_{mj}}
$$
Now for all $m$ let $\hat{a}_m=\hat{T}_{mi}/\hat{T}_{mj}$. Then
\begin{equation}  \label{4-expression}
\frac{F'(\xi)}{|\hat{c}_p-\hat{c}_q|}=\sum_n
\frac{\hat{c}_n-\hat{c}_q}{|\hat{c}_p-\hat{c}_q|} \left(
\frac{\hat{y}_n e^{(\hat{c}_n-\hat{c}_q) \xi} \hat{a}_n
\hat{T}_{nj}}{\sum_m \hat{y}_m e^{(\hat{c}_m-\hat{c}_q) \xi}
\hat{a}_m \hat{T}_{mj}}-\frac{\hat{y}_n e^{(\hat{c}_n-\hat{c}_q)
\xi} \hat{T}_{nj}}{\sum_m \hat{y}_m e^{(\hat{c}_m-\hat{c}_q) \xi}
\hat{T}_{mj}} \right)=\sum_n
\frac{\hat{c}_n-\hat{c}_q}{|\hat{c}_p-\hat{c}_q|} \hat{D}_n.
\end{equation}
where $\hat{D}_n$ denotes the quantity in parentheses in the middle
expression above.

Let $x,y \in W^\circ$ such that for all $k$, $|\hat{x}_k-x_k| \leq
\delta |x_k|$ and  $|\hat{y}_k-y_k| \leq \delta |y_k|$. Let
$a_m=T_{mi}/T_{mj}$,  $c_m=\log x_m/y_m$, and let $D_n$ denote the
unperturbed version of $\hat{D}_n$:
\begin{equation} \label{Dn}
D_n=\frac{y_n e^{(c_n-c_q) \xi} a_n T_{nj}}{\sum_m y_m e^{(c_m-c_q)
\xi} a_m T_{mj}}-\frac{y_n e^{(c_n-c_q) \xi} T_{nj}}{\sum_m y_m
e^{(c_m-c_q) \xi} T_{mj}}.
\end{equation}
By Lemma~\ref{aA}, we have
\begin{equation} \label{posi-nega}
\sum_{n \in \mathcal{T}_0} D_n = \sum_{n \in \mathcal{T}_1} |D_n|
\leq \max_{k, l} \frac{1-\sqrt{a_k/a_l}}{1+\sqrt{a_k/a_l}} \leq
\tau(T),
\end{equation}
where $\mathcal{T}_0=\{n: D_n \ge 0\}$ and $\mathcal{T}_1=\{n: D_n <
0\}$.

Now, for some universal constant $K_0$,
\begin{equation}
\label{real-diff} \left| \sum_n
\frac{\hat{c}_n-\hat{c}_q}{|\hat{c}_p-\hat{c}_q|} \hat{D}_n - \sum_n
\frac{\hat{c}_n-\hat{c}_q}{|\hat{c}_p-\hat{c}_q|} D_n \right| <
K_0(\delta +r).
\end{equation}

Applying Lemma~\ref{maxone} twice, we conclude that there exist $n_0
\in \mathcal{T}_0, n_1 \in \mathcal{T}_1$ such that
$$
\left|\sum_n
\frac{\hat{c}_n-\hat{c}_q}{|\hat{c}_p-\hat{c}_q|} D_n \right|
$$
$$
\leq \left|\frac{\hat{c}_{n_0}-\hat{c}_q}{|\hat{c}_p-\hat{c}_q|}
\sum_{n \in \mathcal{T}_0} D_n + \sum_{n \in \mathcal{T}_1}
\frac{\hat{c}_n-\hat{c}_q}{|\hat{c}_p-\hat{c}_q|} D_n \right|
\leq \left|\frac{\hat{c}_{n_0}-\hat{c}_q}{|\hat{c}_p-\hat{c}_q|}
\sum_{n \in \mathcal{T}_0} D_n -
\frac{\hat{c}_{n_1}-\hat{c}_q}{|\hat{c}_p-\hat{c}_q|} \sum_{n \in
\mathcal{T}_1} |D_n| \right|.
$$
Then together with (\ref{Lijbound}), (\ref{4-expression}),  (\ref{real-diff}), (\ref{posi-nega}),
and the fact that
$|\hat{c}_{n_1}-\hat{c}_{n_0}| \leq |\hat{c}_p-\hat{c}_q|$, we
obtain that for sufficiently small $r, \delta >0$,
$|L_{ij}|$ is upper bounded by some $\rho  < 1$,  as desired.
\end{proof}

\begin{rem}
\label{nbhd} One can further choose $r, \delta > 0$ such that when
$\hat{T} \in B_T(r)$,  $f_{\hat{T}}(W_{\mathbb{C}}^{\circ}(\delta))
\subset W_{\mathbb{C}}^{\circ}(\delta)$. Consider a compact subset
$N \subset W^{\circ}$ such that $f_T(W) \subset N$. Let $N(R)$
denote the Euclidean $R$-neighborhood of $N$ in $W_{\mathbb{C}}$.
The proof of Theorem~\ref{ComplexContraction} implies that when $T >
0$ or ($T \geq 0$ and $\sup_{x,y \in N, 0 \leq \xi \leq 1} \sum_{n
\in \tau_0} D_n < 1$ (here $D_n$ is defined in (\ref{Dn}))), there
exist $r, R > 0$ such that when $\hat{T} \in B_T(r)$, $f_{\hat{T}}$
is a contraction mapping on $N(R)$ under the complex Hilbert metric.
\end{rem}

\begin{exmp} \label{2by2}
Consider a $2 \times 2$ strictly positive matrix
$$
T=\left[ \begin{array}{cc}
a&c\\
b&d\\
\end{array} \right].
$$
If we parameterize the interior of the simplex $W^\circ$ by
$(0,\infty)$: $w = (x,y) \mapsto x/y$, then letting $z = x/y$, we
have: $f_T(z)=\frac{a z+b}{c z+d}$; the domain of this mapping naturally extends from
$(0,\infty)$ to the open right half complex plane $H$, and the
complex Hilbert metric becomes simply $d_H(z_1,z_2) =
|\log(z_1/z_2)|$.

One can show that $f_T$ is a contraction on all of
$H$ with contraction coefficient:
$$
\tau(T)=\frac{1-\frac{bc}{ad}}{1+\frac{bc}{ad}}.
$$
(assuming $\det(T) \ge 0$; otherwise, the last expression is replaced
by $\frac{1-\frac{ad}{bc}}{1+\frac{ad}{bc}}$). To see this, for any $z, w \in H$, consider
$$
L=\left| \frac{\log(f_T(z))-\log(f_T(w))}{\log(z)-\log(w)} \right|.
$$
With change of variables $u=\log(z), v=\log(w)$, we have
$$
L= \left| \frac{\log(f_T(e^u))-\log(f_T(e^v))}{u-v} \right| = \left| \int_0^1 e^{v+t(u-v)} \frac{f'_T(e^{v+t(u-v)})}{f_T(e^{v+t(u-v)})} dt \right|,
$$
which implies that
$$
L \leq \sup_{z \in H} \left| \frac{z f_T'(z)}{f_T(z)} \right|.
$$
A simple computation shows that
\begin{equation}
\label{compute} \frac{z f_T'(z)}{f_T(z)} =  \frac{ad -
bc}{acz+(ad+bc)+bd/z}.
\end{equation}

To see that the supremum is
$\frac{1-\frac{bc}{ad}}{1+\frac{bc}{ad}}$, first note that since $ad
-bc \ge 0$ and $a,b,c,d >0$, the absolute value of the quantity on
the right-hand side of (\ref{compute}) is maximized by minimizing
$|acz + bd/z|$; since the only solutions to $acz + bd/z = 0$ are $z
= \pm i\sqrt{bd/ac}$, one sees that the supremum is obtained by
substituting $z = \pm i\sqrt{bd/ac}$ into~(\ref{compute}), and this
shows that the supremum is indeed
$\frac{1-\frac{bc}{ad}}{1+\frac{bc}{ad}}$.

Note that this contraction coefficient on $H$ is strictly larger
(i.e., worse) than the contraction coefficient on $[0, \infty)$:
$\frac{1-\sqrt{\frac{bc}{ad}}}{1+\sqrt{\frac{bc}{ad}}}$.

When
$$\hat{T} =\left[ \begin{array}{cc}
\hat{a}&\hat{c}\\
\hat{b}&\hat{d}\\
\end{array} \right].
$$
is a sufficiently small complex perturbation of $T$, then $f_{\hat
T}(H) \subseteq H$ and one obtains
$$
\tau(\hat{T})=\sup_{z \in H} \left| \frac{z
f_{\hat{T}}'(z)}{f_{\hat{T}}(z)} \right|=\sup_{z \in H} \left|
\frac{\hat{a}\hat{d} -
\hat{b}\hat{c}}{\hat{a}\hat{c}z+(\hat{a}\hat{d}+\hat{b}\hat{c})+
\hat{b}\hat{d}/z} \right|
$$
which will approximate $\frac{1-\frac{bc}{ad}}{1+\frac{bc}{ad}}$,
and so $f_{\hat{T}}$ will still be a contraction on $H$.
\end{exmp}

\begin{rem}
\label{Dubois}

While this paper was nearing completion, we were informed that
alternative complex Hilbert metrics, based on the Poincare metric in
the right-half complex plane, were recently introduced in
Rugh~\cite{Rugh} and Dubois~\cite{Dubois}. Contractiveness with
respect to these metrics is proven in great generality and yields
far-reaching consequences for complex Perron-Frobenius theory. The
proofs of contractiveness in these papers seem rather different from
the calculus approach in our paper.

The complex Hilbert metric, which we call $d_P$, used
in~\cite{Dubois} (see equation (3.23)) is explicit and natural, but
slightly more complicated than our complex Hilbert metric; for $v,w
\in W_{\mathbb{C}}^+$,
\begin{equation}
\label{dp} d_P(w,v) = \log \frac{\max_{i, j} (|\overline{w_i}v_j +
\overline{w_j}v_i| + |w_iv_j -
w_jv_i|)(2\mathcal{R}(\overline{w_i}w_j))^{-1}} {\min_{i, j}
(|\overline{w_i}v_j + \overline{w_j}v_i| - |w_iv_j -
w_jv_i|)(2\mathcal{R}(\overline{w_i}w_j))^{-1}};
\end{equation}
here $\overline{z}$ denotes complex conjugate, $\mathcal{R}(z)$
denotes real part, and $\log$ is the ordinary real logarithm. In the 2-dimensional case, it can be
verified that, if one transforms $w = (w_1,w_2)$ and $v = (v_1,v_2)$
to $z_1 = w_2/w_1$ and $z_2 =w_2/w_1$, then $d_P$ reduces to the
Poincare metric on $H$:
$$
d_P(z_1,z_2) = \log \frac{|z_1 + \bar{z_2}| + |z_1  - z_2|} {|z_1 +
\bar{z_2}| - |z_1  - z_2|}.
$$
Using the infinitesimal form for the Poincare metric (as a
Riemannian metric on $H$), one checks that, in the $2 \times 2$
case, the Lipschitz constant for a complex matrix $\hat{T}$ such
that $f_{\hat T}(H) \subseteq H$ is:
\begin{equation}
\label{inf-Poincare} \sup_{z \in H} \left| \frac{\mathcal{R}(z)
f_{\hat{T}}'(z)}{\mathcal{R}(f_{\hat{T}}(z))} \right|
\end{equation}
in contrast to
\begin{equation}
\label{inf-Hilbert}
 \sup_{z \in H} \left| \frac{z
f_{\hat{T}}'(z)}{f_{\hat{T}}(z)} \right|
\end{equation}
for our complex Hilbert metric (as in Example~\ref{2by2} above).

While we have not analyzed in detail the differences between these
metrics, there are a few things that can be said in the $2 \times 2$
case:
\begin{itemize}
\item
$f_{\hat{T}}$ is a contraction with respect to $d_P$ on $H$ whenever
it maps $H$ into its interior; this follows from standard complex
analysis (section IX.3 of of~\cite{Gamelin}), and
Dubois~\cite{Dubois} proves an analog of this for the metric $d_P$
above (\ref{dp})  in higher dimensions.  However, this does not hold
for $d_H$.
\item When $\hat{T} = T$ is strictly positive, then the contraction
coefficient, with respect to $d_P$, is always at least as good
(i.e., at most) the contraction coefficient with respect to $d_H$.
This can be seen as follows:

First recall that any fractional linear transformation $T$ can be
expressed as the composition of transitions, dilations and
inversions.  In the case where $T$ is strictly positive, the
translations are by positive real numbers and the dilations are by
real numbers; see page 65 of~\cite{Gamelin}.  Using the
infinitesimal forms (\ref{inf-Poincare}, \ref{inf-Hilbert}), our
assertion would follow from:
\begin{equation}
\label{ineq} \left|\frac{\mathcal{R}(z)}{z}\right| \le
\left|\frac{\mathcal{R}(f_T(z))}{f_T(z)}\right|, ~ \mbox{ for all }
z \in H.
\end{equation}
This is true indeed: it is easy to see that in fact we get equality
in (\ref{ineq}) for inversions and dilations by real numbers,  and
we get strict inequality in (\ref{ineq}) for translations by
positive real numbers.

\item When $\hat{T}$ is a complex perturbation of a strictly
positive $T$, then (\ref{ineq}) (with $T$ replaced by $\hat{T}$)
need not hold; in fact, for perturbations $\hat{T}$ of $T$ on the
order of 1\% and $z = x +yi \in H$, with $|y|/x$ on the order of
1\%, the contraction coefficient with respect to $d_H$ may be
slightly smaller than that with respect to $d_P$. The reason is that
in this case, the dilations may be complex (non-real) and for such a
dilation the inequality (\ref{ineq}) may be reversed. Examples of
this can be randomly generated in Matlab.  For example, if
$$
\hat{T}  =\left[ \begin{array}{cc}
0.012890500224 + 0.000128905002i&  0.310402226067 + 0.003104022260i\\
0.779079247486 - 0.007790792474i
& 0.307296084921 - 0.003072960849i\\
\end{array} \right]
$$
and $z = 0.926678310631 - 0.009266783106i $, then the contraction
coefficent of $d_H$ is approximately  0.664396 and that of $d_P$ is
approximately 0.664599. For larger perturbations, the differences in
contraction coefficient can be greater.
The relative strength of contraction of $d_H, d_P$ seems to be
heavily dependent on specific choices of $\hat{T}$ and $z$.
\item
For any point $z$, other than 0, of the imaginary axis, the metric
$d_H$ can be extended to a neighbourhood, with respect to which any
sufficiently small complex perturbation $\hat{T}$ of a strictly
positive matrix acts as a contraction; on the other hand, there is no way to do this with
$d_P$ since it blows up as one approaches the imaginary axis.
\item Also, on a small punctured neighbourhood of $0$, we replace $d_H$ by the metric
$d(z_1,z_2) = |\log(z_1) - \log(z_2)|$, then small complex perturbation $\hat{T}$ of a strictly
positive matrix still acts as a contraction.
\end{itemize}

In the next section, we use $d_H$ for estimates on the domain of
analyticity of entropy rate of a hidden Markov process.
Alternatively, $d_P$ could be used, however it appears to be
computationally easier to use $d_H$ for the estimation.
\end{rem}

\section{Domain of Analyticity of Entropy Rate of Hidden Markov Processes}

\subsection{Background} \label{background}

For $m, n \in \mathbb{Z}$ with $m \leq n$,  we denote a sequence of
symbols $y_m, y_{m+1}, \ldots, y_n$ by $y_m^n$. Consider a
stationary stochastic process $Y$ with a finite set of states
$\mathcal{I}=\{1, 2, \cdots, B\}$ and distribution $p(y_m^n)$.
Denote the conditional distributions by $p(y_{n+1}|y_m^n)$. The
entropy rate of $Y$ is defined as
$$
H(Y)=\lim_{n \to \infty} -E_p(\log(p(y_0|y_{-n}^{-1}))),
$$
where $E_p$ denotes expectation with respect to the distribution
$p$.

Let $Y$ be a stationary first order Markov chain with
$$
\qquad \Delta(i, j)=p(y_1=j | y_0=i).
$$
In this section, we only consider the case when $\Delta$ is strictly positive.

A \emph{hidden Markov process (HMP)} $Z$ is a process of the form $Z=\Phi(Y)$, where $\Phi$ is a
function defined on $\mathcal{I}=\{1, 2, \cdots, B\}$ with values in
$\mathcal{J}=\{1, 2, \cdots, A\}$.

Recall that $W$ is the $B$-dimensional real simplex and $W_{\mathbb{C}}$ is the complex version of $W$.
For $a \in \mathcal{J}$, let
$\mathcal{I}(a)$ denote the set of all indexes $i \in \mathcal{I}$
with $\Phi(i)=a$.
Let
$$
W_a = \{w \in W: w_i=0  \mbox{ whenever } i \not\in \mathcal{I}(a) \}
$$
and
$$
W_{a, {\mathbb{C}}} = \{ w \in W_{\mathbb{C}}:
w_i=0 \mbox{ whenever } i \not\in \mathcal{I}(a)  \}.
$$

Let $\Delta_a$ denote the $B \times B$ matrix such that $\Delta_a(i,
j)=\Delta(i, j)$ for $j \in \mathcal{I}(a)$, and $\Delta_a(i, j)=0$
for $j \notin \mathcal{I}(a)$ (i.e, $\Delta_a$ is formed from
$\Delta$ by ``zeroing out'' the columns corresponding to indices
that are not in   $\mathcal{I}(a)$. the  For $a \in \mathcal{J}$,
define the scalar-valued and vector-valued functions $r_a$ and $f_a$
on $W$ by
$$
r_a(w)= w \Delta_a \mathbf{1},
$$
and
$$
f_a(w)=w \Delta_a/ r_a(w).
$$
Note that $f_a$ defines the action of the matrix $\Delta_a$ on the
simplex $W$. For any fixed $n$ and $z_{-n}^0$ and for $i=-n, -n+1, \cdots$, define
\begin{equation} \label{x-i}
x_i=x_i(z_{-n}^i)=p(y_i=\cdot \;|z_i, z_{i-1}, \cdots, z_{-n}),
\end{equation}
(here $\cdot$ represent the states of the Markov chain $Y$);
then
from Blackwell~\cite{bl57}, we have that $\{x_i\}$ satisfies the random dynamical
iteration
\begin{equation}
\label{iter0} x_{i+1}=f_{z_{i+1}}(x_i),
\end{equation}
starting with
\begin{equation}
\label{init0} x_{-n-1} = p(y_{-n-1}=\cdot\;).
\end{equation}
where $p(y_{-n-1}=\cdot\;)$ is the stationary distribution for the underlying
Markov chain. One checks that $p(z_{i+1}|z_{-n}^{i})$ can be recovered from this dynamical system; more specifically, we have
$$
p(z_{i+1}|z_{-n}^{i})=r_{z_{i+1}}(x_i).
$$

If the entries of $\Delta = \Delta^{\vec{\varepsilon}}$ are
analytically parameterized by a real variable vector
$\vec{\varepsilon} \in \mathbb{R}^k$ ($k$ is a positive integer), then we obtain a family $Z =
Z^{\vec{\varepsilon}}$ and corresponding $\Delta_a =
\Delta_a^{\vec{\varepsilon}}$, $f_a = f_a^{\vec{\varepsilon}}$, etc.

The following result was proven in~\cite{gm05}.
\begin{thm}  \label{main}
Suppose that the entries of $\Delta = \Delta^{\vec{\varepsilon}}$
are analytically parameterized by a real variable vector
$\vec{\varepsilon}$. If at $\vec{\varepsilon}=\vec{\varepsilon}_0$,
$\Delta$ is strictly positive, then $H(Z) =
H(Z^{\vec{\varepsilon}})$ is a real analytic function of
$\vec{\varepsilon}$ at $\vec{\varepsilon}_0$.
\end{thm}

In~\cite{gm05} this result is stated in greater generality, allowing
some entries of $\Delta$ to be zero.  The proof is based on an
analysis of the action of perturbations of  $f_a$ on neighbourhoods
of $\hat{W}_b\stackrel{\triangle}{=}f_b(W)$, with respect to the Euclidean metric. The proof assumes
that each $f_a$ is a contraction on each $\hat{W}_b$. While this need not
hold, one can arrange for this to be true by replacing the original
system with a higher power system: namely, one replaces the original
alphabet $\mathcal{J}$ with $\mathcal{J}^n$ for some $n$ and
replaces the mappings $\{f_a: a \in \mathcal{J}\}$ with $\{f_{a_0}
\circ f_{a_1} \circ \cdots \circ f_{a_{n-1}}: a_0 a_1 \ldots a_{n-1} \in
\mathcal{J}^n\}$. The existence of such an $n$ follows from a) the
equivalence of the (real) Hilbert metric and the Euclidean metric on
each $\hat{W}_b$ (Proposition 2.1 of ~\cite{gm05}) and b) the
contractiveness of each $f_a$ with respect to the (real) Hilbert
metric.  However, in the course of this replacement, one easily
loses track of the domain of analyticity.

When at $\vec{\varepsilon}=\vec{\varepsilon}_0$, $\Delta$ is
strictly positive, an alternative is to directly use a complex
Hilbert metric, as follows. For each $a \in \mathcal{J}$, we can
define a complex Hilbert metric $d_{a, H}$ on $W^\circ_{a,
\mathbb{C}}$ as follows: for $w, v \in W^\circ_{a, \mathbb{C}}$:
\begin{equation}
d_{a, H}(w, v)=d_H(w_{\mathcal{I}(a)},v_{\mathcal{I}(a)})=\max_{i, j \in \mathcal{I}(a) } \left| \log \left(
\frac{w_i/w_j}{v_i/v_j} \right) \right|.
\end{equation}
Theorem~\ref{ComplexContraction} implies that for each $a,b \in
\mathcal{J}$, sufficiently small perturbations of $f_a$ are
contractions on sufficiently small complex neighborhoods of
$\hat{W}_b$ in $W_{b, \mathbb{C}}$; see Remark~\ref{nbhd} (note that
while $\Delta_a$ is not strictly positive, $f_a$ maps into $W_a$ and
so as a mapping from $W_b$ to $W_a$ it can be regarded as the
induced mapping of a strictly positive matrix). For complex
$\vec{\varepsilon}$ close to $\vec{\varepsilon}_0$, $f_a =
f_a^{\vec{\varepsilon}}$ is sufficiently close to
$f_a^{\vec{\varepsilon_0}}$ to guarantee that
$f_a^{\vec{\varepsilon}}$ is a contraction.

Let $\Omega_{a,H}(R)$ denote the neighborhood of diameter $R$, measured in the complex Hilbert metric,
of $\hat{W}_a$ in $W_{a, \mathbb{C}}$. Let $B_{\vec{\eps}_0}(r)$ denote the complex $r$-neighborhood of
$\vec{\eps}_0$ in $\mathbb{C}^k$.

Following the proof of Theorem~\ref{main}
(especially pages 5254-5255 of~\cite{gm05}), one obtains a lower bound $r>0$ on
the domain of analyticity if there exists $R > 0$ and  $0 < \rho
<1$ satisfying the following conditions:
\begin{enumerate}
\item For any $a, z \in \mathcal{A}$ and any $\vec{\eps} \in B_{\vec{\eps}_0}(r)$,
$f_z^{\vec{\eps}}$ is a contraction, with respect to the complex
Hilbert metric, on $\Omega_{a,H}(R)$:
$$
\sup_{x \neq y \in \Omega_{a,H}(R)} \left| \frac{d_{z, H}
(f_z^{\vec{\eps}}(x), f_z^{\vec{\eps}}(y))}{d_{a, H}(x, y)} \right| \leq \rho < 1.
$$
\item
for any
$\vec{\eps} \in B_{\vec{\eps}_0}(r)$, any $x \in \cup_a \hat{W}_a$ and any $z \in \mathcal{A}$,
$$
d_{z, H}(f_z^{\vec{\eps}}(x), f_z^{\vec{\eps}_0}(x)) \leq R(1-\rho),
$$
and
$$
d_{z, H}(f_z^{\vec{\varepsilon}}(\pi(\eps)), f_z^{\vec{\varepsilon}_0}(\pi(\eps_0))) \leq R(1-\rho),
$$
(where $\pi(\eps)$ denotes the stationary vector for the Markov chain defined
by $\Delta^{\vec{\varepsilon}}$).
\item
For any $x \in \Omega_{{a,H}}(R)$ and $\vec{\eps} \in B_{\vec{\eps}_0}(r)$,
$$
\sum_a |r^{\vec{\eps}}_a(x)| \leq 1/\rho.
$$
\end{enumerate}

The existence of $r, R, \rho$ follows from Theorem~\ref{main}. In
fact, we can choose $\rho$ to be any positive number such that $\max_{a
\in \mathcal{A}} \tau (\Delta_a) < \rho < 1$, and small $r, R$ to
satisfy condition 1, then smaller $r, R$, if necessary, to further
satisfy conditions 2 and 3.

Let $\Omega_{a,E}(R)$ denote the neighborhood of diameter
$R$, measured in the Euclidean metric, of $\hat{W}_a$ in $W_{a,
{\mathbb{C}}}$. To facilitate the computation, at the expense of
obtaining a smaller lower bound, it may be easier to use
$\Omega_{a,E}(R)$ instead of $\Omega_{a,H}(R)$; then, the conditions
above are replaced with the following conditions:
\begin{itemize}
\item[(1')] Condition 1 above with $\Omega_{a,H}(R)$ replaced by $\Omega_{a,E}(R)$
(the map $f_z^{\vec{\eps}}$ is still required to be a contraction
under the complex {\em Hilbert} metric).
\item[(2')]
Condition 2 above with $R$ on the right hand side of the
inequalities replaced by $R/K$, where $K=\sup_{x \neq y \in
\Omega_{a,E}(R), a} \left| \frac{d_{a,E}(x, y)}{d_{a, H}(x, y)}
\right|$; note that for $R$ sufficiently small, $0 < K < \infty$
since $d_{a, H}$ and $d_{a, E}$ are equivalent metrics (this in turn
follows from the fact that the Euclidean metric and (real) Hilbert metric are equivalent
on any compact subset of the interior of the real simplex).
\item[(3')] Condition 3 above with $\Omega_{a,H}(R)$ replaced by $\Omega_{a,E}(R)$
\end{itemize}

\subsection{Example for Domain of Analyticity}
\label{example}

In the following, we consider hidden Markov processes obtained by
passing binary Markov chains through binary symmetric channels with
crossover probability $\eps$. Suppose that the Markov chain is
defined by a $2 \times 2$ stochastic matrix $\Pi = [\pi_{ij}]$. From
now through the end of this section, we {\bf assume:}
\begin{itemize}
\item  $\det(\Pi) > 0$ -- and --
\item all $\pi_{ij} > 0$ -- and --
\item $0 < \varepsilon < 1/2$.
\end{itemize}
We remark that the condition $\det(\Pi) > 0$ is purely for
convenience.

Strictly speaking, the underlying Markov process of the resulting hidden Markov process is given by a 4-state
matrix (the states are the ordered pairs of a state of $\Pi$ and a
noise state (0 for ``noise off'' and 1 for ``noise on''); see page 5255 of~\cite{gm05}). However,
the information contained in each $f_a$ can be reduced to an
equivalent map induced by a $2 \times 2$ matrix and then reduced to
an equivalent function of a single variable as in
Example~\ref{2by2}. We describe this as follows.

Let $a_i=p(z_1^i, y_i=0)$ and $b_i=p(z_1^i, y_i=1)$. The pair
$(a_i,b_i)$ satisfies the following dynamical system:
$$
(a_i, b_i)=(a_{i-1}, b_{i-1}) \left[ \begin{array}{cc}
p_E(z_i) \pi_{00} & p_E(z_i) \pi_{10}\\
p_E(\bar{z}_i) \pi_{01} & p_E(\bar{z}_i) \pi_{11}\\
\end{array} \right].
$$
where $p_E(0) = \eps$ and $p_E(1) = 1 - \eps$.

Similar to Example~\ref{2by2}, let $x_i=a_i/b_i$, we have a dynamical system with just one
variable:
$$
x_{i+1}=f^\eps_{z_{i+1}}(x_i),
$$
where
$$
f^\eps_{z}(x)=\frac{p_E(z)}{p_E(\bar{z})} \frac{\pi_{00}
x+\pi_{10}}{\pi_{01} x+\pi_{11}}, \qquad  z=0, 1
$$
starting with
\begin{equation} \label{x0}
x_0=\pi_{10}/\pi_{01},
\end{equation}
which comes from the stationary vector of $\Pi$.

It can be shown that
$$
p^\eps(z_i=0|z_1^{i-1})=r^\eps_0(x_{i-1}), \qquad
p^\eps(z_i=1|z_1^{i-1})=r^\eps_1(x_{i-1}),
$$
where
\begin{equation} \label{r0}
r^\eps_0(x)= \frac{((1-\varepsilon) \pi_{00}+\varepsilon \pi_{01})
x+((1-\varepsilon) \pi_{10}+ \varepsilon \pi_{11})}{x+1},
\end{equation}
and
\begin{equation}  \label{r1}
r^\eps_1(x)=\frac{(\varepsilon \pi_{00}+(1-\varepsilon) \pi_{01})
x+(\varepsilon \pi_{10}+ (1-\varepsilon) \pi_{11})}{x+1}.
\end{equation}


Now let $\Omega(R)$ denote the complex $R$-neighborhood (in Euclidean metric) of the
interval
$$
S=[S_1, S_2]=\left[ \frac{\eps_0 \pi_{10}}{(1-\eps_0) \pi_{11}},
\frac{(1-\eps_0) \pi_{00}}{\eps_0 \pi_{01}} \right],
$$
this interval is the union of $f^{\eps_0}_0([0,\infty])$ and
$f^{\eps_0}_1([0,\infty])$; again let $B_{\vec{\eps}_0}(r)$ denote the complex
$r$-neighborhood of a given cross-over probability $\eps_0>0$.

The sufficient conditions (1'), (2') and (3') in
section~\ref{background} are guaranteed by the following: there exist
$R
>0, r>0, 0< \rho <1$ such that
\begin{enumerate}
\item[(1'')] For any $z$, $f_z^{\eps}(x)$ is a contraction on $\Omega(R)$ under complex Hilbert metric,
$$
\sup_{x \neq y \in \Omega(R)} \left| \frac{\log f_z^{\eps}(x)-\log f_z^{\eps}(y)}{\log x-\log y} \right| \leq
\rho < 1.
$$
Note that here
$$
\log f_z^{\eps}(x)-\log f_z^{\eps}(y)=\log \frac{\pi_{00}
x+\pi_{10}}{\pi_{01} x+\pi_{11}} -\log \frac{\pi_{00}
y+\pi_{10}}{\pi_{01} y+\pi_{11}}.
$$
\item[(2'')] For any $\eps \in B_{\vec{\eps}_0}(r)$, any $x \in S$ and any $z$,
$$
|\log f_z^{\eps}(x)-\log f_z^{\eps_0}(x)| \leq (R/K)(1-\rho),
$$
where
$$
K=\sup_{x \neq y \in \Omega(R)} \left| \frac{x-y}{\log x -\log y}
\right|=\sup_{x \in \Omega(R)} |x|=S_2+R.
$$
(note that here the second condition in (2') is vacuous since by (\ref{x0}) $x_0$ does not depend on $\eps$)
\item[(3'')] For any $x \in \Omega(R)$ and $\eps \in B_{\vec{\eps}_0}(r)$,
$$
|r^\eps_0(x)|+|r^\eps_1(x)| \leq 1/\rho.
$$
\end{enumerate}

By considering extreme cases, the above conditions can be further
relaxed to:
\begin{enumerate}
\item[(1''')]
$$
0 < \frac{\pi_{00} \pi_{11}-\pi_{10} \pi_{01}}{\pi_{01} \pi_{00}
(S_1-R)+\pi_{01} \pi_{10}+\pi_{11} \pi_{00} + \pi_{11} \pi_{10}/(S_2
+R)} \leq \rho.
$$
(here we applied the mean value theorem to give an upper bound on $|\log((\pi_{00}
x+\pi_{10})/(\pi_{01} x+\pi_{11})) -\log((\pi_{00}
y+\pi_{10})/(\pi_{01} y+\pi_{11}))|$)

\item[(2''')]
$$
0 < \frac{r}{\eps_0-r}+\frac{r}{1-\eps_0-r} \leq (R/(S_2+R))(1-\rho).
$$
(here we applied the mean value theorem to give an upper bound on $|\log((1-\eps)/\eps)-\log((1-\eps_0)/\eps_0)|$)

\item[(3''')]
$$
0 < \frac{((1-\eps_0+r) \pi_{00}+(\eps_0+r) \pi_{01})
(S_2+R)+((1-\eps_0+r) \pi_{10}+ (\eps_0+r) \pi_{11})}{S_1-R+1}
$$
$$
+\frac{((\eps_0+r) \pi_{00}+(1-\eps_0+r) \pi_{01})
(S_2+R)+((\eps_0+r) \pi_{10}+ (1-\eps_0+r) \pi_{11})}{S_1-R+1} \leq
1/\rho.
$$
\end{enumerate}

In other words, choose $r, R$ and $\rho$ to satisfy the conditions
(1'''), (2''') and (3'''). Then the entropy rate is an analytic
function of $\eps$ on $|\eps-\eps_0| < r$.

\begin{figure}
\begin{center}
\includegraphics[totalheight=3in]{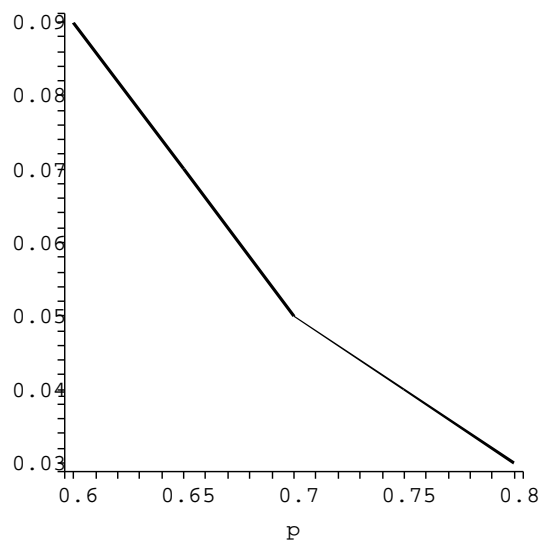}
\caption{lower bound on radius of convergence as a function of $p$}
\label{LowerBound}
\end{center}
\end{figure}

Consider the symmetric case: $\pi_{00}=\pi_{11}=p$ and
$\pi_{01}=\pi_{10}=1-p$. We plot lower bounds on radius of
convergence of $H(Z)$ (as a function of $\eps$ at $\eps_0=0.4$)
against $p$ in Figure~\ref{LowerBound}. For a
fixed $p$, the lower bound is obtained by randomly generating many
$3$-tuples $(r, R, \rho)$ and taking the maximal $r$ from the
$3$-tuples which satisfy the conditions.

\bigskip

{\bf Acknowledgements:} We thank Albert Chau for helpful discussions
on Riemannian metrics in $H$.


\begin{thebibliography}{10}

\bibitem{bl57}
D.~Blackwell.
\newblock The entropy of functions of finite-state {M}arkov chains.
\newblock {\em Trans. First Prague Conf. Information Thoery, Statistical
  Decision Functions, Random Processes}, pages 13--20, 1957.

\bibitem{Dubois}
L.~Dubois.
\newblock Projective metrics and contraction principles for complex cones.
\newblock {\em http://arxiv.org/abs/0811.2930}.

\bibitem{Gamelin}
T.~Gamelin.
\newblock Complex Analysis
\newblock Springer, 2001.

\bibitem{Holliday}
Holliday, T., Goldsmith, A., and Glynn, P.
\newblock Capacity of Finite
State Channels Based on Lyapunov Exponents of Random Matrices
\newblock{\em IEEE Transactions on Information Theory}, Volume 52, Issue 8,
August 2006, pages 3509 - 3532

\bibitem{gm05}
G.~Han and B.~Marcus.
\newblock Analyticity of entropy rate of hidden {M}arkov chains.
\newblock {\em IEEE Transactions on Information Theory},
Volume 52, Issue 12, December, 2006, pages: 5251-5266.

\bibitem{hm06b}
G.~Han and B.~Marcus.
\newblock Derivatives of Entropy Rate in
Special Familes of Hidden Markov Chains.
\newblock {\em IEEE Transactions on Information Theory}, Volume 53, Issue 7,
July 2007, Page(s):2642 - 2652.

\bibitem{or03}
E.~Ordentlich and T.~Weissman.
\newblock On the optimality of symbol by symbol filtering and denoising.
\newblock {\em IEEE Transactions on Information Theory}, Volume 52, Issue 1, Jan. 2006, Page(s):19 -
40.

\bibitem{Peres1} Y. Peres.
\newblock
Analytic dependence of Lyapunov exponents on transition
probabilities,
\newblock
Spriner Lecture Notes in Mathematics, Lyapunov's exponents,
Proceedings of a Workshop, volume 1486, Springer Verlag, 1990.

\bibitem{Peres2}
Y. Peres. \newblock Domains of analytic continuation for the top
Lyapunov exponent. \newblock Ann. Inst. H. Poincare Probab.
Statist., 28(1):131--148, 1992.

\bibitem{ro03}
Rolando Cavazos-Cadena.
\newblock An alternative derivation of Birkhoff's formula
for the contraction coefficient of a positive matrix.
\newblock Linear Algebra and its Applications, 375, (2003), 291-297.

\bibitem{Rugh}
H.~H.~Rugh.
\newblock Cones and gauges in complex spaces : Spectral gaps and complex Perron-Frobenius theory.
\newblock {\em http://arxiv.org/abs/math/0610354}.

\bibitem{se80}
E.~Seneta.
\newblock {\em Springer Series in Statistics}.
\newblock Non-negative Matrices and Markov Chains. Springer-Verlag, New York
Heidelberg Berlin, 1980.

\bibitem{ta02}
J.~L.~Taylor.
\newblock {\em Several complex variables with connections to algebraic geometry
and Lie groups}.
\newblock American Mathematical Society, Providence, R.I., 2002.

\end{thebibliography}
\end{document}